\documentclass{article}
\usepackage{amsmath, amssymb, graphics, setspace}

\newcommand{\mathsym}[1]{{}}
\newcommand{\unicode}[1]{{}}

\begin{document}

\title{Un exemple de fonction int{\' e}grable sur \([-1,1]\) mais pas sur \([0,1]\)}
\author{\text{Fran{\c c}ois Gu{\' e}nard}}
\date{}
\maketitle

\textbf{\textit{Abstract $-$ In courses on integration theory, Chasles property is usually considered as elementary and so natural that the is sometimes left to the reader.  When the functions take their values in finite dimensional spaces, the property is always verified, but it no more true in infinite dimensional spaces.  We first give an easy-to-understand example of a function \textit{f}  from   \([-1, 1]\)  into the space of polynomial functions from   \([0, 1]\)  to $\mathbb{R}$ which is integrable on  \([-1, 1]\)  but not on  \([0, 1]\). We also provide a way of representing graphically such a function which explains what means the integral of a function with values in an infinite dimensional space.  Then we show that Chasles'property is true if and only if the space in which the functions to integrate take their values is a complete space. }}\\

Dans les livres traitant de l{'}int{\' e}gration de fonctions d{'}une variable r{\' e}elle {\` a} valeurs dans un espace de dimension infinie \textit{
E}, il est d{'}usage de supposer que l{'}espace d{'}arriv{\' e}e de la fonction est complet. { }Cela permet d{'}assurer la convergence des sommes
de Riemann, et d{'}avoir le th{\' e}or{\` e}me de Chasles qui affirme, pour une fonction \(f : [-1, 1] \rightarrow  E\), les assertions suivantes
sont {\' e}quivalentes :\\
(i)$\quad $\textit{ f} est int{\' e}grable sur \([-1, 1]\) ;\\
(ii)$\quad $il existe un r{\' e}el \(c \in  ]-1, 1[\) tel que \textit{ f} soit int{\' e}grable sur \([-1, c]\) et sur \([c, 1]\); \\
(iii) $\quad $pour tout r{\' e}el { }\(c \in  ]-1, 1[\), \textit{ f} est int{\' e}grable sur \([-1, c]\) et sur \([c, 1]\).

La d{\' e}monstration du th{\' e}or{\` e}me de Chasles repose sur l{'}utilisation du crit{\` e}re de Cauchy pour la convergence des sommes de Riemann,
selon un type de subdivision d{\' e}pendant de la th{\' e}orie de l{'}int{\' e}grale consid{\' e}r{\' e}e : Riemann, Lebesgue (McShane), Denjoy,
Henstock. { }

Dans cette note, on va tout d{'}abord donner un exemple explicite de fonction d{\' e}finie sur \([-1, 1]\), {\` a} valeurs dans un espace vectoriel
norm{\' e} \textit{ E}, int{\' e}grable sur \([-1, 1]\), mais non int{\' e}grable sur \([-1, 0]\), ni sur \([0, 1]\). On donnera {\' e}galement une
illustration graphique de cette fonction. { }On montrera dans une deuxi{\` e}me partie que ceci est g{\' e}n{\' e}ral : Le th{\' e}or{\` e}me de
Chasles est valide si et seulement si l{'}espace d{'}arriv{\' e}e est complet.

\section{Construction explicite d{'}un contre-exemple au th{\' e}or{\` e}me de Chasles}

\subsection{Introduction}

L{'}exemple que l{'}on va consid{\' e}rer ici est une fonction impaire de \([-1, 1]\) dans l{'}espace \(E = \mathbb{R}[x, [0, 1]]\) des fonctions
polynomiales de \([0, 1]\) dans $\mathbb{R}$ muni de la norme de la convergence uniforme. { }Comme cet espace vectoriel admet une base alg{\' e}brique
d{\' e}nombrable, il r{\' e}sulte du th{\' e}or{\` e}me de Baire qu{'}il n{'}est complet pour aucune norme. { }On va {\' e}tablir le r{\' e}sultat
pour les sommes de Riemann sur des subdivisions r{\' e}guli{\` e}res. { }Le compl{\' e}t{\' e} de \textit{ E} {\' e}tant l{'}espace \(F = \mathcal{C}([0,
1], \mathbb{R})\) des fonctions continues de \([0, 1]\) dans $\mathbb{R}$ muni de la norme de la convergence uniforme, la convergence sur \([-1,
0]\) et sur \([0, 1]\) des sommes de Riemann r{\' e}guli{\` e}res de \textit{ f} vers deux fonctions \textit{ g} et \textit{ h} appartenant {\` a}
\(F \backslash  E\) montrera l{'}int{\' e}grabilit{\' e} au sens de Riemann de \textit{ f }en tant que fonction de \([-1, 1]\) dans \textit{ F}.
Compte tenu des liens entre les diff{\' e}rentes th{\' e}ories (\(\text{}^1\)), cela montrera {\' e}galement l{'}int{\' e}grabilit{\' e} de \textit{
f} aux sens de Lebesgue (via l{'}int{\' e}grale de McShane), de Denjoy et de Henstock. { }L{'}int{\' e}grabilit{\' e} de \textit{ f} dans l{'}un
de ces sens sur \([-1, 0]\) et sur \([0, 1]\) avec une int{\' e}grale appartenant {\` a} \(F \backslash  E\) montrera {\' e}galement qu{'}en tant
que fonction de \([-1, 1]\) dans \textit{ E}, cette fonction n{'}est int{\' e}grable en aucun de ces sens sur \([-1, 0]\) ni sur \([0, 1]\). { }
En revanche, l{'}int{\' e}grabilit{\' e} de \textit{ f} sur \([0, 1]\) au sens de Riemann en tant que fonction {\` a} valeurs dans \textit{ E} montrera
son int{\' e}grabilit{\' e} sur \([0, 1]\) aux sens de Lebesgue, Denjoy et Henstock. { }

\subsubsection{Comment repr{\' e}senter concr{\` e}tement l{'}int{\' e}grale d{'}une fonction {\` a} valeurs dans un espace de dimension infinie}

Afin d{'}illustrer graphiquement la fonction \textit{ f}, on va d{\' e}finir dans l{'}ordinateur une fonction de deux variables, $\Phi $, { }de \([-1,
1] \times  [0, 1]\) dans $\mathbb{R}$,\((x , t )|\rightarrow  (f(x))(t)\). { }Cette fonction de deux variables pourra {\^ e}tre repr{\' e}sent{\'
e}e par une surface, ce qui permettra de repr{\' e}senter concr{\` e}tement la fonction \textit{ f} d{\' e}finie sur \([-1, 1]\), {\` a} valeurs
dans un espace de dimension infinie. Une section de la surface par le plan vertical d{'}{\' e}quation \(x = x_0\) montrera le graphe de \(f\left(x_0\right)\).
{ }Pour illustrer les sommes de Riemann de cette fonction, on r{\' e}alisera une animation montrant l{'}{\' e}volution de ces sommes de Riemann.
{ }Pour chacune de ces sommes, on mettra en couleur les graphes des \(f\left(x_i\right)\), o{\` u} les \(x_i\) sont les {\' e}tiquettes de la somme
de Riemann. { }La somme de Riemann elle-m{\^ e}me, une fonction de \textit{ E}, sera repr{\' e}sent{\' e}e sur une face visible de la {``}bo{\^ \i}te{''}
d{\' e}limitant le graphique. { }

\subsection{D{\' e}finition de la fonction}

On va commencer par construire \textit{ f} par morceaux sur \([0, 1]\), sur les intervalles de la forme \(I_n = \left[\frac{1}{2^{n+1}}, \frac{1}{2^n}\right]\),
\(n \geq 0\). { }L{'}id{\' e}e est de reconstituer la fonction exponentielle {\` a} partir de son d{\' e}veloppement en s{\' e}rie. L{'}int{\' e}grale
sur \(I_n\) de \textit{ f} devra valoir \(g_n : [0, 1] \rightarrow  \mathbb{R}, x |\rightarrow  \frac{x^n}{n!}\). { }Pour cela, on va prendre \textit{
f} affine par morceaux sur \(I_n:\)Si \(e_n\)est la fonction de \([-1, 1]\) dans $\mathbb{R}$ d{\' e}finie par \(e_n : x |\rightarrow  x^n\), et
si $\theta $ est la fonction nulle de \([-1, 1]\) dans $\mathbb{R}$, \(f\) est la fonction impaire de \([-1, 1]\) dans \textit{ E} d{\' e}finie par

\(\begin{array}{c}
 
\begin{array}{c}
  f : x |\rightarrow  \left\{
\begin{array}{cc}
 \frac{2^{2n+4}}{n!}\left(x - \frac{1}{2^{n+1}}\right)e_n  & \text{si} x \in  \left[\frac{1}{2^{n+1}}, \frac{3}{2^{n+2}}\right] \\
 \frac{2^{2n+4}}{n!}\left(\frac{1}{2^n}- x\right)e_n & \text{si} x \in  \left[\frac{3}{2^{n+2}}, \frac{1}{2^n}\right] \\
 \theta  & \text{si} x =0
\end{array}
\right.
\end{array}

\end{array}\)

\subsection{Repr{\' e}sentation}

Afin d{'}illustrer graphiquement la fonction \textit{ f}, on a repr{\' e}sent{\' e} la fonction de deux variables, $\Phi $, { }de \([-1, 1] \times
 [0, 1]\) dans $\mathbb{R}$,\((x , t )|\rightarrow  (f(x))(t)\). { }La surface associ{\' e}e {\` a} cette fonction de deux variables permet de repr{\'
e}senter concr{\` e}tement la fonction \textit{ f} d{\' e}finie sur \([-1, 1]\), {\` a} valeurs dans un espace de dimension infinie. Une section
de la surface par le plan vertical d{'}{\' e}quation \(x = x_0\) montre le graphe de \(f\left(x_0\right)\). { }Cette repr{\' e}sentation n{'}est
pas reproduite dans cette version de la note.

\begin{doublespace}
\noindent\(\pmb{\text{}}\)
\end{doublespace}

\subsection{D{\' e}monstration}

La fonction \textit{ f} est continue et affine par morceaux sur ]0, 1] qui est un ouvert de \([0, 1]\). { }Elle est donc continue sur \(]0, 1]\).
{ }En outre, comme \(\frac{2^{2n+4}}{n!}\) tend vers 0 lorsque \(n \rightarrow \infty\), \textit{ f} tend vers la fonction nulle $\theta $ en 0.
Comme \(\theta \in  E\), \textit{ f} est continue sur \([0, 1]\). { }Par parit{\' e}, \textit{ f} est aussi continue sur \([-1, 1]\).

Comme \textit{ f} est affine par morceaux, son int{\' e}grale sur chaque intervalle o{\` u} elle est affine se calcule pratiquement comme l{'}int{\'
e}grale d{'}une fonction {\` a} valeurs r{\' e}elles. { }Voici pourquoi.

\pmb{ Lemme }$--$ Soit \((E, \|.\|)\) un espace vectoriel norm{\' e} r{\' e}el non r{\' e}duit {\` a} $\{$0$\}$, \textit{ e} un vecteur non nul
de \textit{ E}, et \textit{ f} la fonction affine de \([0, 1]\) dans \textit{ E}, \(x |\rightarrow  x e\). Alors \textit{ f} est int{\' e}grable
sur \([0, 1]\) aux sens de Riemann, Henstock, Lebesgue et Denjoy, et \(\int _{[0,1]}f=\frac{1}{2}e\).

\pmb{ D{\' e}monstration du lemme $--$ }Compte tenu des inclusions entre th{\' e}ories de l{'}int{\' e}grale, il suffit d{'}{\' e}tablir que \textit{
f} est int{\' e}grable au sens de Riemann, donc pour des subdivisions r{\' e}guli{\` e}res. { }Soient \(n\in \mathbb{N}\) et \(\sigma (f, n) = \frac{1}{n}\sum
_{k=1}^n f\left(\frac{k}{n}\right)= \frac{n+1}{2n}e\). { }D{'}apr{\` e}s les propri{\' e}t{\' e}s de la norme, on a \(\|\sigma (f, n)\| = \frac{n+1}{2n}\|e\|\).
{ } Comme la droite vectorielle \(\mathbb{R} e\) est compl{\` e}te, \(\sigma (f, n)\) converge dans \(E\) vers \(\frac{1}{2}e\). { }$\blacksquare
$

Le lemme s{'}{\' e}tend {\' e}videmment aux fonctions affines quelconques. { }En l{'}appliquant {\` a} chacun des intervalles \(\left[\frac{1}{2^{n+1}},
\frac{3}{2^{n+2}}\right]\), on trouve que\\

\(\begin{array}{c}
 
\begin{array}{c}
 \int _{\left[\frac{1}{2^{n+1}}, \frac{1}{2^n}\right]}f= 
\end{array}
\int _{\left[\frac{1}{2^{n+1}}, \frac{3}{2^{n+2}}\right]}\frac{2^{2n+4}}{n!}\left(x - \frac{1}{2^{n+1}}\right)e_n+\int _{\left[\frac{3}{2^{n+2}},
\frac{1}{2^n}\right]}\frac{2^{2n+4}}{n!}\left(\frac{1}{2^n}- x\right)e_n
\\
= \frac{2^{2n+4}}{n!}\left(\int _{\left[\frac{1}{2^{n+1}}, \frac{3}{2^{n+2}}\right]}\left(x
- \frac{1}{2^{n+1}}\right)+\int _{\left[\frac{3}{2^{n+2}}, \frac{1}{2^n}\right]}\left(\frac{1}{2^n}- x\right)\right)e_n
\\
=\frac{2^{2n+4}}{n!}\left(\frac{1}{2^{2n+5}}+\frac{1}{2^{2n+5}}\right)
e_n = \frac{1}{n!}e_n
\end{array}\)\\

Ainsi, sur tout segment \(J \subset  ]0, 1]\), \(f\) est int{\' e}grable dans \textit{ E} puisque son int{\' e}grale est une fonction polyn{\^ o}me.
{ }En outre, \(\int _{\frac{1}{2^{n+1}}}^1f=\sum _{k=0}^n \frac{1}{k!}e_k\). { }La restriction {\` a} \([0, 1]\) de la fonction exponentielle appartient
{\` a} \(F\), mais pas {\` a} \(E\). { }Par suite, comme \(\sum _{k=0}^n \frac{1}{k!}e_k\) converge uniform{\' e}ment vers la restriction {\` a}
\([0, 1]\) de la fonction exponentielle, la fonction \textit{ f} est int{\' e}grable sur \([0, 1]\) en tant que fonction de \([0, 1]\) dans \(F\),
mais pas en tant que fonction de \([0, 1]\) dans \(E\). { }La fonction \textit{ f} {\' e}tant impaire, son int{\' e}grale dans \textit{ F} est la
fonction nulle, $\theta $. { }Comme \(\theta \in  E\), la fonction \textit{ f} est donc, en tant que fonction de \([-1, 1]\) dans \(E\) int{\' e}grable
sur \([-1, 1]\), mais n{'}est int{\' e}grable ni sur \([-1, 1]\) ni sur \([0, 1]\). { }

On notera au passage que \(f\) est un exemple de fonction continue sur un segment non int{\' e}grable sur ce segment : l{'}implication {``}\(f\)
continue sur un segment{''} $\Rightarrow $ {``} \(f\) est int{\' e}grable sur ce segment{''} est vraie lorsque l{'}espace d{'}arriv{\' e}e est complet,
fausse sinon.

\section{Cas g{\' e}n{\' e}ral}

Soit \(E\) un espace vectoriel norm{\' e} r{\' e}el non complet, \(F\) son compl{\' e}t{\' e} pour la norme consid{\' e}r{\' e}e. { }Nous allons
construire une fonction continue \(f\) de \([-1, 1]\) dans \(E\) int{\' e}grable sur \([-1, 1]\), mais non int{\' e}grable, ni sur \([-1, 0]\), ni
sur \([0, 1]\). { }

Puisque \(E\) n{'}est pas complet, il existe une suite de Cauchy \(\left(x_n\right)\) d{'}{\' e}l{\' e}ments de \(E\) qui converge vers un {\' e}l{\'
e}ment \(y\) de \(F \backslash  E\). { }Quitte {\` a} extraire des sous-suites, on peut supposer que, pour tout \(n \in  \mathbb{N}^*\), on a \(\left\|y-x_n\right\|
\leq  \frac{1}{2^{2n+3}}\). { }On a alors , pour tout entier { }\(n \geq  2\), \(\left\|x_n-x_{n-1}\right\| \leq \frac{1}{2^{2n}}\). { }Cela {\'
e}tant, on va construire une fonction \textit{ f} impaire, continue, affine par morceaux, et telle que pout tout \(n \in \mathbb{N}^*\), { }\(\int
_{\frac{1}{2^n}}^1f=x_n\). { }Soit \textit{ f} la fonction impaire de \([-1, 1]\) dans \(E\) d{\' e}finie par\\

\(\begin{array}{c}
 
\begin{array}{c}
 f(x) = \begin{cases}
 16(1-x)x_1 & \text{si} x \in  \left[\frac{3}{4}, 1\right] \\
 16\left(x-\frac{1}{2}\right)x_1 & \text{si} x \in  \left[\frac{1}{2}, \frac{3}{4}\right] \\
 2^{2n+2}\left(\frac{4}{2^{n+1}}-x\right)\left(x_n- x_{n-1}\right) & \text{si} x \in  \left[\frac{3}{2^{n+1}}, \frac{4}{2^{n+1}}\right] \\
 2^{2n+2}\left(x-\frac{1}{2^n}\right)\left(x_n- x_{n-1}\right) & \text{si} x \in  \left[\frac{1}{2^n}, \frac{3}{2^{n+1}}\right] \\
 0 & \text{si} x = 0
\end{cases}

\end{array}

\end{array}\)

La fonction \(f\) est continue sur l{'}ouvert \(]0, 1]\) de \([0, 1]\) parce qu{'}elle est affine par morceaux et continue sur tout segment inclus
dans cet intervalle. { }La condition { }\(\left\|x_n-x_{n-1}\right\| \leq \frac{1}{2^{2n}}\) assure la continuit{\' e} de \(f\) en 0. { }Ainsi,\(f\)
est-elle continue sur \([0, 1]\). Gr{\^ a}ce au lemme de la section { }pr{\' e}c{\' e}dente, on montre que \(f\) est int{\' e}grable sur tout segment
\(J \subset \text{  }]0, 1]\). { }En outre, { } \(\int _{\frac{1}{2^n}}^1f=x_n\). { }Ainsi, dans \(F\), \(f\) est int{\' e}grable sur \([0, 1]\),
et \(\int _{(F)[0,1]}f= y\). { }Mais par d{\' e}finition, \(y \notin  E\), et \(f\), en tant que fonction de \([-1, 1]\) dans \(E\) n{'}est pas int{\'
e}grable sur \([0, 1]\). { }Par parit{\' e}, elle n{'}est pas non plus int{\' e}grable sur \([-1, 0]\). { }On a donc {\' e}tabli que

\textbf{Le th{\' e}or{\` e}me de Chasles pour les int{\' e}grales est vrai si et seulement si  
l{'}espace d{'}arriv{\' e}e est complet.}

\subsection{R{\' e}f{\' e}rences}

(\(\text{}^1\)) G. A. Gordon : The integrals of Lebesgue, Denjoy, Perron and Henstock. { }Graduate Studies in Math. Vol 4, Amer Math Soc. 1994.

\begin{doublespace}
\noindent\(\pmb{\text{}}\)
\end{doublespace}

Fran{\c c}ois Gu{\' e}nard\\
Laboratoire de Math{\' e}matiques d{'}Orsay\\
UMR 8628 du CNRS\\
Universit{\' e} Paris-Sud\\
91405 Orsay Cedex, France\\
francois.guenard@math.u-psud.fr

\end{document}